\documentclass[11pt]{article}

\usepackage{latexsym,color,graphicx,amsmath,amssymb}

\def\reals{{\mathbb R}}
\def\R{{\mathbb R}}
\def\eps{{\varepsilon}}

\newtheorem{theorem}{Theorem}[section]
\newtheorem{lemma}[theorem]{Lemma}
\newtheorem{corollary}[theorem]{Corollary}

\def\cplx{\mathbb{C}}
\def\P{\mathbb{P}}

\def\C{{\cal C}}
\def\F{{\cal F}}

\title{On rich points and incidences with restricted sets of lines in 3-space\thanks{%
Partially supported by ISF Grant 260/18, by grant 1367/2016
from the German-Israeli Science Foundation (GIF), and by
Blavatnik Research Fund in Computer Science at Tel Aviv University.
A preliminary version appeared in the Proceedins of the
37th Symposium on Computational Geometry (2021), 56:1--56:14, and in arXiv:2012.11913.}}

\author{
Micha Sharir\thanks{%
School of Computer Science, Tel Aviv University, Tel Aviv Israel;
{\tt michas@tauex.tau.ac.il}.}
\and
Noam Solomon\thanks{%
School of Computer Science, Tel Aviv University, Tel Aviv Israel;
{\tt noam.solom@gmail.com}.}
}

\begin{document}
 
\maketitle

\begin{abstract}
Let $L$ be a set of $n$ lines in $\R^3$ that is contained, when represented
as points in the four-dimensional Pl\"ucker space of lines in $\R^3$, 
in an irreducible variety $T$ of constant degree which is \emph{non-degenerate} 
with respect to $L$ (see below). We show:

\medskip
\noindent{\bf (1)}
If $T$ is two-dimensional, the number of $r$-rich points (points incident to at least $r$ lines of $L$)
is $O(n^{4/3+\eps}/r^2)$, for $r \ge 3$ and for any $\eps>0$, and, if at most $n^{1/3}$ lines of $L$ 
lie on any common regulus, there are at most $O(n^{4/3+\eps})$ $2$-rich points.
For $r$ larger than some sufficiently large constant, the number of $r$-rich points is also $O(n/r)$.

As an application, we deduce (with an $\eps$-loss in the exponent) the bound obtained by 
Pach and de Zeeuw~\cite{PdZ} on the number of distinct distances determined by $n$ points 
on an irreducible algebraic curve of constant degree in the plane that is not a line nor a circle.

\medskip
\noindent{\bf (2)}
If $T$ is two-dimensional, the number of incidences between $L$ and a set of $m$ points 
in $\R^3$ is $O(m+n)$. 

\medskip
\noindent{\bf (3)}
If $T$ is three-dimensional and nonlinear, the number of incidences between $L$ and a set of $m$ points 
in $\R^3$ is $O\left(m^{3/5}n^{3/5} + (m^{11/15}n^{2/5} + m^{1/3}n^{2/3})s^{1/3} + m + n \right)$,
provided that no plane contains more than $s$ of the points. When $s = O(\min\{n^{3/5}/m^{2/5}, m^{1/2}\})$, 
the bound becomes $O(m^{3/5}n^{3/5}+m+n)$.

As an application, we prove that the number of incidences between $m$ points and $n$ lines 
in $\R^4$ contained in a quadratic hypersurface (which does not contain a hyperplane) 
is $O(m^{3/5}n^{3/5} + m + n)$.

The proofs use, in addition to various tools from algebraic geometry, recent bounds on 
the number of incidences between points and algebraic curves in the plane.
\end{abstract}

\section{Introduction} \label{sec:intro}

\paragraph{The setup: Incidences between a set of points and a restricted set of lines in $\R^3$.}
Let $P$ be a set of $m$ points and $L$ a set of $n$ lines in $\R^3$. 
We consider the problem of obtaining sharp incidence bounds between the points of $P$ and the lines of $L$, 
when the lines of $L$, considered as points in the four-dimensional Pl\"ucker space of lines in $\R^3$,
are restricted to lie on a two- or three-dimensional constant-degree algebraic variety $T$.
The topic of incidences between points and lines is a fundamental topic in incidence geometry, 
significantly boosted since Guth and Katz's seminal work~\cite{GK2} on point-line incidences in $\R^3$.
Instead of asking for a bound on the number of incidences between points and lines, we can ask,
for each $r\ge 3$, for a bound on the number of \emph{$r$-rich points} in a set of lines, 
which are the points that are incident to at least $r$ of the lines. As it turns out, the 
two questions are equivalent. The related, and finer problem of bounding the number of $2$-rich points, 
determined by a set of $n$ lines in $\R^3$, studied in \cite{GK2}, is also discussed in this paper,
under the restricted setup considered here.
Building on recent works of Sharir and Zahl~\cite{SZ} and Zahl~\cite{Za17}, we are able to improve 
Guth and Katz's point-line incidence bounds when the lines in $L$ are restricted to lie on a two- or 
three-dimensional variety $T$ in the Pl\"ucker space. 

\paragraph{Background: Points and curves, the planar case.}
The study of incidences between points and curves has a rich history,
starting with the simplest instance of points and lines in the
plane, where we have (see also~\cite{CEGSW,Sze}):
\begin{theorem}[Szemer\'edi and Trotter \cite{ST}] \label{th:ST}
The maximum number of incidences between $m$ points and $n$ lines in the plane
is $\Theta(m^{2/3}n^{2/3}+m+n)$.
\end{theorem}
In fact, an equivalent formulation of Szemer\'edi-Trotter theorem asserts that, 
given $n$ lines in the plane, the number of points that are incident 
to at least $r$ of the lines, for any parameter $2 \le r \le n$, which we call
\emph{$r$-rich points} and denote the set of these points by $P_{\ge r}(L)$, satisfies
\begin {equation}
\label{eq:szt}
|P_{\ge r}(L)| = O\left( \frac{n^2}{r^3} + \frac{n}{r} \right).
\end {equation}
Still in the plane, Pach and Sharir~\cite{PS} extended this bound to incidence
bounds between points and curves with $k$ \emph{degrees of freedom}, namely,
for each set of $k$ distinct points,
there are only $\mu=O(1)$ curves that pass through all of them, and
each pair of curves intersect in at most $\mu$ points; $\mu$ is
called the \emph{multiplicity} (of the degrees of freedom).
Here is their result, tailored to the case of algebraic curves.
\begin{theorem}[Pach and Sharir \cite{PS}] \label{th:PS}
Let $P$ be a set of $m$ points in $\R^2$ and let $\C$ be a set
of $n$ bounded-degree algebraic curves in $\R^2$ with $k$
degrees of freedom and with multiplicity $\mu$. Then
(where the constant of proportionality depends on $k$ and $\mu$)
\[
I(P,\C) = O\left(m^{\frac{k}{2k-1}}n^{\frac{2k-2}{2k-1}}+m+n\right) .
\]
\end{theorem}

\medskip
Except for the case $k=2$ (lines have two degrees of freedom), the upper
bound is not known, and is strongly suspected to be too large 
(see~\cite{ANPPSS, ArS, MT} for an improvement for the case of circles and similar curves).

Recently, Sharir and Zahl~\cite{SZ} have considered general families
of constant-degree algebraic curves in the plane that belong to an
\emph{$s$-dimensional family of curves}. This means that each curve in 
such a family can be represented by a constant number of real parameters, 
so that, in this parametric space, the points representing the curves lie 
in an $s$-dimensional algebraic variety $\F$ of some constant degree
(the so-called ``complexity'' of $\F$). See~\cite{SZ} for details.
\begin{theorem}[Sharir and Zahl~\cite{SZ}]\label{incPtCu}
Let $\C$ be a set of $n$ algebraic plane curves that belong to
an $s$-dimensional family $\F$ of curves of maximum constant degree
$E$, no two of which share a common irreducible component, and let
$P$ be a set of $m$ points in the plane. Then, for any $\eps>0$,
the number $I(P,\C)$ of incidences between the points of
$P$ and the curves of $\C$ satisfies
\begin{equation*} \label{newIncBd}
 I(P,\C) = O\Big(m^{\frac{2s}{5s-4}} n^{\frac{5s-6}{5s-4}+\eps} + m^{2/3}n^{2/3} + m + n\Big) ,
\end{equation*}
where the constant of proportionality depends on $\eps$, $s$, $E$,
and the complexity of the family $\F$.
\end{theorem}
Except for the factor $O(n^\eps)$, this is a significant improvement
over the bound in Theorem~\ref{th:PS} (for $s\ge 3$), in cases where
the assumptions in Theorem~\ref{incPtCu} imply (as they often do)
that $\C$ has $k=s$ degrees of freedom. Examples where $k=s$ are 
abundant in the plane. For example, lines have $k=s=2$ (two points 
determine a line, and two real parameters specify a line. For circles 
we have $k=s=3$, for unit circles we have $k=s=2$, and for general
conic sections we have $k=s=5$. 

\paragraph{Incidences with lines in three dimensions.}

The groundbreaking work of Guth and Katz~\cite{GK2} implies\footnote{%
  This bound is not explicitly stated in \cite{GK2}, but it readily
  follows from the analysis given there, and by now it is generally
  attributed to that work.}
the sharper bound $O(m^{1/2}n^{3/4} + m^{2/3}n^{1/3}q^{1/3} + m + n)$ 
on the number of incidences between $m$ points and $n$ lines in $\R^3$,
provided that no plane contains more than $q$ of the given lines.
We use the following variant of this result.
\begin{theorem} \label{th:ss4d}
Let $P$ be a set of $m$ points in $\R^3$, and let $L$ be a set of $n$ lines 
in $\R^3$, so that no 2-flat contains more than $s$ points of $P$. Then
\[
I(P,L) = O(m^{1/2}n^{3/4} + m^{1/3}n^{2/3}s^{1/3} + m + n).
\] 
\end{theorem}
Note that this symmetric assumption (no more than $s$ points, instead
of no more than $q$ lines, on a plane) causes the roles of $m$ and $n$ 
to switch in the second term, from their roles in the Guth--Katz bound.

\noindent{\bf Sketch of proof.}
The proof is very similar to that in \cite{GK2}. Briefly, it is
based on constructing a \emph{partitioning polynomial} $f$ of
some suitable degree. It then bounds the number of incidences
within the cells of the partition (connected components of
$\R^3\setminus Z(f)$, where $Z(f)$ is the zero set of $f$), 
and the number of incidences on $Z(f)$. The latter task is
performed separately for each irreducible component of $Z(f)$.
We follow the proof in \cite{GK2} verbatim, except for the analysis
of incidences on the \emph{planar} irreducible components of $Z(f)$.
For that step, both arguments use the Szemer\'edi-Trotter bound
for point-line incidences in the plane (Theorem~\ref{th:ST}),
for each of these planes, except that the analysis in \cite{GK2}
exploits the assumption that no plane contains more than $q$ of the lines, whereas
we exploit the symmetric assumption that no plane contains more than $s$ of the points.
Reasoning as in \cite{GK2}, but swapping the roles of points and lines, we obtain
the modified second term in the bound in the theorem. See Section~\ref{sec:inc:3d} 
for a similar reasoning for the restricted cases considered here.
$\Box$

Plugging the bound of Theorem \ref{th:ss4d} into the proof of \cite[Theorem 1.3(a)]{SS4dv}, we get 
\begin{theorem} \label{th:ss4dv}
Let $P$ be a set of $m$ points in $\R^3$, and let $L$ be a set of $n$ 
lines in $\R^d$, for $d \ge 3$, so that all
the points and lines lie in a common two-dimensional algebraic variety $V$ of degree $D$ that
does not contain any 2-flat, and so that no 2-flat contains more than $s$ points of $P$. Then
\[
I(P,L) = O(m^{1/2}n^{1/2}D^{1/2} + m^{1/3}D^{4/3}s^{1/3} + m + n).
\] 
\end{theorem}

Guth and Katz's work has lead to many recent works on incidences 
between points and lines or other curves in three and higher dimensions;
see~\cite{CPS,GZ,SSS,SSZ,SS4d,SS4dv} for a sample.

Of particular significance is the recent work of Guth and
Zahl~\cite{GZ} on the number of 2-rich points in a collection of
algebraic curves of constant degree, namely, points incident to at 
least two of the given curves. For the case of lines, Guth and 
Katz~\cite{GK2} have shown that the number of such points is 
$O(n^{3/2})$, when no plane or regulus contains more than 
$O(n^{1/2})$ lines. Guth and Zahl obtain the same asymptotic 
bound for general algebraic curves, under analogous (but
stricter) restrictive assumptions (concerning surfaces that
are doubly ruled by the given family of curves).

The new bounds that we will derive require the extension to three dimensions 
of the notions of having $k$ degrees of freedom and of being an 
$s$-dimensional family of curves. The definitions of these concepts, 
as given above for the planar case, extend, more or less verbatim, to 
three (or higher) dimensions, but, even in typical situations, these 
two concepts do not coincide anymore. For example, lines in three 
dimensions have two degrees of freedom, but they form a $4$-dimensional 
family of curves (this is the number of parameters needed to
specify a line in $\R^3$). 

\paragraph{Our results.}
We obtain improved incidence bounds when the lines of $L$, as points
in Pl\"ucker space, lie on a two- or three-dimensional variety $T$.
When $T$ is two-dimensional and non-planar, the number of $r$-rich points 
is $O(n^{4/3+\eps}/r^2)$, for $r \ge 3$ and for any $\eps>0$, and, 
if at most $n^{1/3}$ lines of $L$ lie on any common regulus, there 
are at most $O(n^{4/3+\eps})$ $2$-rich points. For $r$ larger than 
some sufficiently large constant, the number of $r$-rich points is 
also $O(n/r)$, which is a better bound when $r = O(n^{1/3})$. 
These bounds improve significantly, for the restricted context at hand, the bound
$O(n^{3/2}/r^2)$ due to Guth and Katz~\cite{GK2} (which holds
when no plane or regulus contains more than $O(n^{1/2})$ lines).
Moreover, the number of incidences between $L$ and a set of $m$ 
points in $\R^3$ is $O(m+n)$, again a significant improvement, in our context,
over the previous bound in \cite{GK2}.

As an application, we show that the number of distinct distances 
determined by $n$ points on an irreducible algebraic curve of 
constant degree in the plane that is not a line nor a circle,
is $\Omega(n^{4/3-\eps})$, for any $\eps>0$, which is
(with an $\eps$-loss in the exponent) the bound obtained by 
Pach and de Zeeuw~\cite{PdZ}.

If $T$ is three-dimensional and nonlinear, the number of incidences between $L$ and a set of $m$ points 
in $\R^3$ is $O\left(m^{3/5}n^{3/5} + (m^{11/15}n^{2/5} + m^{1/3}n^{2/3})s^{1/3} + m + n \right)$,
provided that no plane contains more than $s$ of the points. When $s = O(\min\{n^{3/5}/m^{2/5}, m^{1/2}\})$, 
the bound becomes $O(m^{3/5}n^{3/5}+m+n)$.

An interesting novel feature of our result is that, like Theorem~\ref{th:ss4d}, 
it is obtained under an assumption that restricts the number of \emph{points} 
that can lie on a common plane (instead of restricting the number of coplanar 
lines in the previous studies). Very few earlier works have used this kind of 
restriction; see Elekes et al.~\cite{EKS} for one of the few exceptions.

Similar bounds have recently been obtained by the authors for other special 
cases of the incidence problem~\cite{SSZl,SZ}, using related but different approaches.

As an application, we prove that the number of incidences between $m$ points and $n$ lines 
in $\R^4$ contained in a quadratic hypersurface (which does not contain a hyperplane) 
is $O(m^{3/5}n^{3/5} + m + n)$.

All our bounds are significant improvements, under the restricted scenarios 
assumed in this work, over the standard incidence bounds in
three dimensions, and shed, as we believe, new light on the structure of point-line 
incidences in three dimensions.

As is standard in the `modern' study of incidence geometry, the analysis
is based on the \emph{polynomial partitioning} technique (see~\cite{Guth,GK2} for details),
combined with a variety of tools from algebraic geometry.

We remark that, wherever needed in the analysis, we switch to 
the (projective 3-space over) the complex field, which simplifies it 
and lets us use numerous tools from algebraic geometry, available
in this setting. The passage from the complex projective 
setup back to the real affine one is easy---the former is a 
generalization of the latter. The real affine setup is needed only for 
constructing a polynomial partitioning, which is meaningless over
$\cplx$. Once we are, say, within the zero set $Z(f)$ of
the partitioning polynomial $f$, we can switch to the
complex projective setup, and reap the benefits just noted.

\section{Rich points determined by a two-dimensional family of lines}
\label{rich-2d}

As already said, we parameterize lines in three dimensions 
by their \emph{Pl\"ucker coordinates}, as follows (see, e.g.,
Griffiths and Harris~\cite[Section 1.5]{GrHa}). For a pair of
distinct points $x,y \in \P^3$, given in projective coordinates 
as $x=(x_0,x_1,x_2,x_3)$ and $y=(y_0,y_1,y_2,y_3)$, let 
$\ell_{x,y}$ denote the (unique) line in $\P^3$ incident 
to both $x$ and $y$. The Pl\"ucker coordinates of
$\ell_{x,y}$ are given in projective coordinates in $\P^5$ as
$(\pi_{0,1},\pi_{0,2},\pi_{0,3},\pi_{2,3},\pi_{3,1},\pi_{1,2})$,
where $\pi_{i,j}=x_iy_j-x_jy_i$. Under this parameterization, the set
of lines in $\P^3$ corresponds bijectively to the set of points in
$\P^5$ lying on the \emph{Klein quadric} $Q$ given by the quadratic equation 
\begin{equation} \label{eq:klein}
\pi_{0,1}\pi_{2,3}+\pi_{0,2}\pi_{3,1}+\pi_{0,3}\pi_{1,2} = 0
\end{equation} 
(which is indeed always satisfied by the Pl\"ucker coordinates of a line).

As we will be working both in the primal three-dimensional space and the dual
five-dimensional projective Pl\"ucker space, or rather in the Klein quadric $Q$
contained in that space, we will denote the Pl\"ucker point image of a line $\ell$ 
by $\tilde{\ell}$ throughout the paper. For a set $L$ of lines in $\reals^3$ or
in $\cplx^3$, we denote\footnote{%
  As mentioned in the introduction, one can freely pass between a real algebraic
  variety, or other constructs, and its complexification, so, for convenience,
  we will mostly argue over the complex domain.}
$\tilde{L} = \{\tilde{\ell} \mid \ell\in L\}$.

Given a surface $V$ in $\P^3$, the set of lines fully
contained in $V$, represented by their Pl\"ucker coordinates in
$\P^5$, is a subvariety of the Klein quadric $Q$, which is denoted by $F(V)$,
and is called the \emph{Fano variety} of $V$; see
Harris~\cite[Lecture 6, page 63]{Har} for details, and \cite[Example 6.19]{Har}
for an illustration, and for a proof that $F(V)$ is indeed a variety.

Let $H$ denote a plane in $\cplx^3$, and let $H^* = \{\tilde{\ell} \mid \ell\subset H\}$ 
be the 2-flat in Pl\"ucker space consisting of the points that represent the 
lines fully contained in $H$ (see Rudnev~\cite{Rud} for why $H^*$ is
indeed a 2-flat and for more details).

For a set $L$ of lines, we put 
\[
\mathcal H(L) = \{H_{\ell, \ell'}^* \mid \ell, \ell' \in L \text{ and are coplanar}\},
\] 
where for coplanar lines $\ell, \ell'$, $H_{\ell, \ell'}$ is the (unique) 2-flat 
containing $\ell$ and $\ell'$.

In this paper, we study incidences between a set of points $P\subset \R^3$, 
and a set of lines $L$ in $\R^3$, whose Pl\"ucker images lie on some irreducible 
algebraic subvariety $T$ of the Klein quadric $Q$, which is of constant degree, 
and which has dimension either $2$ or $3$.

In this section we restrict ourselves to the case where $\dim(T) = 2$.
For a set $L$ of lines, we say that the (two-dimensional) variety $T$ 
is \emph{non-degenerate} with respect to $L$ if 
\begin{description}
\item[(i)]
$T$ is irreducible of constant degree, 
\item[(ii)]
$T$ is not a $2$-flat, and 
\item[(iii)]
the intersection of $T$ with each 2-flat $H^* \in \mathcal H(L)$ consists 
of $O(1)$ points (representing lines). 
\end{description}
Note that condition (iii) is what one would expect
to hold in a generic situation in a four-dimensional space. The simpler case 
where $T$ is a $2$-flat will not be considered in this work. The reason is 
that when many of the input lines lie in a plane $h$ in the primal 3-space, 
the dual of $h$ in the Pl\"ucker space is a 2-flat, which could be $T$ itself. 
In this case the best upper bound that one could get for the number of $r$-rich 
points would be the worst-case bound $O(n^2/r^3 + n/r)$ of
Szemer\'edi and Trotter~\cite{ST}.


\begin{theorem}
\label{th:tworich}
Let $L$ be a set of $n$ lines in $\R^3$, such that, in Pl\"ucker space, 
$L$ is a subset of some two-dimensional variety $T$ that is non-degenerate 
with respect to $L$. 

\medskip
\noindent
(a) The number of $r$-rich points determined by $L$ 
is $|P_{\ge r}(L)| = O(n^{4/3 +\eps}/r^2)$, for any $\eps>0$ and $r \ge 3$. 

\medskip
\noindent
(b) If, in addition, the number of lines of $L$ contained in any common regulus
is at most $n^{1/3}$ then the number of $2$-rich points determined by $L$
is $|P_{\ge 2}(L)| = O(n^{4/3 + \eps})$, for any $\eps>0$.
\end{theorem}

\noindent{\bf Proof.}
First here is a high-level overview of the proof. After a pruning step, we
may assume that the set $\gamma_\ell$, for a line $\ell\in L$, of the 
Pl\"ucker images of the lines $\ell$ that are coplanar with $\ell$ and 
$\tilde{\ell}\in T$, is a one-dimensional curve in $T$.
An $r$-rich point generates $\Omega(r^2)$ incidences between the Pl\"ucker
points in $\tilde{L}$ and the curves $\gamma_\ell$, so it suffices to
bound the number of such incidences. There are two kinds of curves, those
that represent the lines in one ruling of some regulus, and those that do not.
For $r$-rich points, with $r\ge 3$, only the latter kind of curves matter,
and a suitable application of Theorem~\ref{incPtCu} allows us to obtain an upper bound
for the number of such incidences. For $2$-rich points (part (b) of the theorem),
the regulus-curves also play a part, and the analysis is complicated because these
curves do not have to be distinct. Still, the assumptions in (b) allow us to handle
this case and get the desired bound.

The following notation will be useful later on in the paper. 
For each line $\ell \in Q$, define the variety $S_\ell$ (in Pl\"ucker space) to be
\[
S_\ell = \{\tilde{\ell}' \in Q \mid \text{ $\ell$, $\ell'$ are coplanar}\} .
\]
If the Pl\"ucker coordinates of $\tilde{\ell}$ are 
$(\pi_{0,1},\pi_{0,2},\pi_{0,3},\pi_{2,3},\pi_{3,1},\pi_{1,2})$, then
\begin{multline*}
S_{\ell} = \{(\pi'_{0,1},\pi'_{0,2},\pi'_{0,3},\pi'_{2,3},\pi'_{3,1},\pi'_{1,2})\in Q \mid \\
\pi_{0,1}\pi'_{2,3}+\pi_{0,2}\pi'_{3,1}+\pi_{0,3}\pi'_{1,2} + \pi'_{0,1}\pi_{2,3}+\pi'_{0,2}\pi_{3,1}+\pi'_{0,3}\pi_{1,2} = 0\} .
\end{multline*}
In particular, Equation~(\ref{eq:klein}) implies that $\tilde{\ell} \in S_{\ell}$.
We see that, for every line $\ell$, the variety $S_\ell$ is the intersection 
of $Q$ with a hyperplane, so it is a three-dimensional quadratic surface 
contained in $Q$ and containing $\tilde{\ell}$.
We say that a line $\ell$ is \emph{exceptional} with respect to $T$ 
if $T \subset S_\ell$. We say that a point $p \in \R^3$ is \emph{exceptional} 
with respect to $T$ if the set of the Pl\"ucker images of the lines incident 
to $p$ in 3-space, which we denote by $S_p$ and which is known to be a 2-plane 
in $Q$ (see the proof of Lemma~\ref{lem:noexcep} below), is equal to $T$. 

\begin{lemma} \label{lem:noexcep}
(i) If $T$ is non-degenerate, there are no exceptional points with respect to $T$. \\
(ii) Even when $T$ is degenerate, there can be at most one exceptional point.
\end{lemma}

\noindent{\bf Proof.}
We recall, e.g., from Rudnev~\cite{Rud}, that a 2-flat contained in $Q$ 
parameterizes either the set of lines in $\cplx^3$ that are incident to some point,
or the set of lines contained in a plane in $\cplx^3$. The set of lines that are 
incident to an exceptional point is thus a 2-flat, which, as $T$ is an irreducible 
two-dimensional surface, must be equal to $T$, contradicting the assumption 
that $T$ is non-degenerate, so (i) follows. For the proof of (ii), we observe 
that if $p, p'$ are two exceptional points, then $S_p, S_{p'} \subset T$ and 
$S_p \ne S_{p'}$, contradicting the assumption that $T$ is two-dimensional and irreducible. 
$\Box$

\begin{lemma}
\label{le:plre}
There are at most two exceptional lines with respect to $T$.
\end{lemma}
\noindent{\bf Proof.}
Assume to the contrary that there are three exceptional lines. Assume first that 
two of these lines are coplanar, and denote them by $\ell_1$ and $\ell_2$. 
Then $T \subseteq S_{\ell_1} \cap S_{\ell_2}$, i.e., $T$ is contained in the set 
of the Pl\"ucker images of the lines intersecting both $\ell_1$ and $\ell_2$. 
If $\ell_1$ and $\ell_2$ do not intersect one another, then 
$S_{\ell_1} \cap S_{\ell_2} = H_{\ell_1, \ell_2}^*$. 
Otherwise, letting $p = \ell_1 \cap \ell_2$, we have
$S_{\ell_1} \cap S_{\ell_2} = H_{\ell_1, \ell_2}^* \cup S_p$, 
i.e., it is a union of two $2$-flats. In both cases,
$T$ is a $2$-flat, contradicting our assumption.

We may thus assume there are (at least) three lines $\ell_1, \ell_2$ and $\ell_3$ 
that are pairwise skew, such that $T \subseteq S_{\ell_1} \cap S_{\ell_2} \cap S_{\ell_3}$. 
As is well known (see, e.g.,~\cite[Theorem 16.4]{FT} and~\cite[Lemma 2.2]{SS4dv}), 
the Pl\"ucker images of lines that intersect $r\ge 3$ pairwise-skew lines 
$\ell_1,\ldots, \ell_r$ belong to one ruling of a regulus, and the Pl\"ucker 
images of $\ell_1,\ldots,\ell_r$ belong to the other ruling of this regulus. 
That is, $S_{\ell_1} \cap S_{\ell_2} \cap S_{\ell_3}$ is one ruling of 
the regulus generated by the lines intersecting $\ell_1, \ell_2$ and $\ell_3$, 
which is a quadratic curve in the Pl\"ucker space, contradicting the fact that 
$T$ is two-dimensional. It is a curve because the ruling defines a one-parameter
family of lines, which yields a curve in Pl\"ucker space; the fact that this curve is
quadratic, an easy consequence of the fact that the regulus spans a quadratic surface,
is not important for the argument.
$\Box$

\medskip

We prune away, as we may, the exceptional point, if such a point exists 
(in the degenerate situation), and the (at most) two exceptional lines, 
thereby losing at most $2(n-1)+1 < 2n$ $2$-rich points. 

For each of the (remaining) lines $\ell \in L$, the intersection $S_\ell \cap T$ 
is a curve contained in $T$ (possibly also containing a discrete finite subset), 
which we denote by $\gamma_\ell$. Define 
\begin{equation} \label{eq:c}
\C = \{\gamma_\ell \mid \ell \text{ is not exceptional}\}.
\end{equation} 
We have the following simple observation, whose trivial proof is omitted.
\begin{lemma}
\label{le:trivial}
Let $p$ be an $r$-rich point, with $r \ge 2$, and denote the lines
of $L$ incident to $p$ as $\ell_1,\ldots,\ell_s$, for some $s\ge r$. 
Then, for each pair of indices $1\le i\ne j\le s$, $\tilde{\ell}_i$
is incident to $\gamma_{\ell_j}$, 
and for every incidence $\tilde{\ell}'\in \gamma_{\ell}$, for a pair of lines
$\ell$, $\ell'$, there is at most one point $p\in \cplx^3$ that induces 
it, in the sense that both $\ell$ and $\ell'$ are incident to $p$.
\end{lemma}

The lemma asserts that each $r$-rich point contributes at least $r(r-1)$ 
incidences between the Pl\"ucker points of the lines of $L$ and the curves
$\gamma_\ell$ of $\C$ (as curves in $Q$). Hence, to bound the number of 
$r$-rich points, it suffices to bound the number of incidences between 
the Pl\"ucker images of the lines in $L$ and the curves of $\C$ 
(and then divide the bound by $r(r-1)$).
For any curve $\gamma_\ell$, its corresponding discrete subset of $O(1)$ points
contributes only $O(1)$ incidences, for a total of $O(n)$ incidences. We
may therefore ignore all these discrete subsets.

The notion of \emph{dimensionality} for families of curves (see the 
definition preceding Theorem~\ref{incPtCu}) easily extends in a natural 
way to collections $\C$ of higher-dimensional algebraic varieties.
\begin{lemma}
\label{le:twodim}
The family $\C$ defined in (\ref{eq:c}) is two-dimensional.
\end{lemma}

\noindent{\bf Proof.}
Each curve $\gamma_\ell$ of $\C$ can be parameterized by the 
parameters of the corresponding line $\ell$, and the Pl\"ucker images
of these lines lie in the two-dimensional variety $T$, so it takes only
two real parameters to specify such a line.
$\Box$

\medskip

When analyzing incidences between the Pl\"ucker points in
$\tilde{L}$ and the curves $\gamma_j$, as in Lemma~\ref{le:trivial}, 
some care has to be exercised, to handle situations in which many of the 
curves $\gamma_\ell$ share a common irreducible component (or even coincide).

\paragraph{$(\ge 3)$-rich points.}
Assume that $\ell_1,\ldots, \ell_\xi \in L$ are such that 
$\gamma_{\ell_1}, \ldots, \gamma_{\ell_\xi}$ all share a common curve, for some $\xi\ge 3$. 
If some pair of lines $\ell_i, \ell_j$ are coplanar, we write $H_{i,j}$ for the (unique) 
plane $H_{\ell_i,\ell_j}$ containing both $\ell_i$ and $\ell_j$. 
As in the proof of Lemma~\ref{le:plre},
(i) if $\ell_i$ and $\ell_j$ are parallel then $S_{\ell_i}\cap S_{\ell_j} = H_{i,j}^*$, 
and (ii) if $\ell_i$ and $\ell_j$ intersect in a point $p$ then
$S_{\ell_i}\cap S_{\ell_j} = H_{i,j}^*\cup S_p$, so $S_{\ell_i}\cap S_{\ell_j}$ 
is either a $2$-flat or the union of two 2-flats (in $Q$). Therefore, 
\[
\gamma_{\ell_i} \cap \gamma_{\ell_j} = S_{\ell_i}\cap S_{\ell_j} \cap T = 
H_{i,j}^* \cap T \quad\text{or}\quad (H_{i,j}^* \cup S_p) \cap T,
\] 
and the right hand sides of these equations are the intersection of one or two 
$2$-flats with $T$. Since $T$ is assumed to be non-degenerate, it 
follows that $\gamma_{\ell_i} \cap \gamma_{\ell_j}$ is a finite set of points, and 
thus $\gamma_{\ell_i}$ and $\gamma_{\ell_j}$ cannot intersect in a common curve.
We can thus assume that $\ell_1,\ldots, \ell_\xi$ are pairwise skew (and 
that $\gamma_{\ell_1},\ldots, \gamma_{\ell_\xi}$ intersect in a common curve).

\begin{lemma}
\label{le:com} 
Assume that the arc (in Pl\"ucker space) $\gamma:=\bigcap_{i=1}^\xi \gamma_{\ell_i}$ 
is nonempty (and is not a finite set), where $\ell_1,\ldots,\ell_\xi$ are 
$\xi$ pairwise-skew lines, with $\xi \ge 3$.
Then $\gamma$ parameterizes one ruling of a regulus, and, for each line $\ell\in T$ 
such that $\tilde{\ell}\in\gamma$, $\ell$ intersects $\ell_1,\ldots, \ell_\xi$, 
and thus $\tilde{\ell}$ lies in the curve that represents the other ruling 
of the same regulus.
\end{lemma}
\noindent{\bf Proof.}
The proof is similar to the proof of Lemma~\ref{le:plre}. 
The intersection $\bigcap_{i=1}^\xi S_{\ell_i}$ consists of the Pl\"ucker points
that represent the lines that intersect the $\xi\ge 3$ pairwise-skew lines 
$\ell_1,\ldots,\ell_\xi$. Thus, as already noted (see~\cite{FT}), all these lines
belong to one ruling of a regulus, and the Pl\"ucker points of $\ell_1,\ldots,\ell_\xi$ 
belong to the other ruling of this regulus. Therefore, $\gamma$ parameterizes one 
ruling of a regulus, and $\ell_1,\ldots,\ell_\xi$ belong to the other ruling of this 
regulus, as asserted.
$\Box$

\medskip

Partition the set of irreducible components of the curves $\gamma_{\ell}$, 
over all lines $\ell \in L$ that are not exceptional, into two subsets $\C_0$ and $\C_1$, 
where $\C_0$ (resp., $\C_1$) contains all the components that do not (resp., do)
parameterize one ruling of some regulus. Since $\deg(\gamma_{\ell}) \le \deg(T) = O(1)$, 
for each $\ell\in L$, it follows that $|\C_0| = |\C_1| = O(n)$. 
We partition the set of incidences into incidences between the Pl\"ucker points 
in $\tilde{L}$ and the curves in $\C_0$, and incidences with
the curves in $\C_1$.
We remind the reader that at this stage we are only concerned with incidences 
induced by a concurrence of at least $r\ge 3$ lines of $L$ at some ($r$-rich) point $p$. 

By Lemma~\ref{le:trivial}, any $r$-rich point $p$, for $r \ge 3$,
corresponds to incidences between the Pl\"ucker points $\tilde{\ell}$ that
represent lines $\ell$ of $L$ that are incident to $p$ and the (at least three)
curves $\gamma_{\ell'}$ that are associated with these lines, and any such 
incidence can arise for at most one point $p$.  
One possibility is that the Pl\"ucker point $\tilde{\ell}$, for a line $\ell$ 
incident to $p$, is incident to a common component of at least three of these 
curves, call them $\gamma_{\ell_1^p}, \gamma_{\ell_2^p}, \gamma_{\ell_3^p}$. 
However, the analysis preceding Lemma~\ref{le:com} 
implies that $\ell_1^p$, $\ell_2^p$, $\ell_3^p$ are pairwise skew, 
which is impossible as they are all incident to $p$.
Hence an incidence between the Pl\"ucker image of a line and a common component 
of at least three curves $\gamma_{\ell_i}$ does not generate any
$r$-rich points, for $r \ge 3$, and, by construction, curves in $\C_1$ also do not 
generate any $r$-rich points, for $r \ge 3$. We may therefore assume that 
every curve in $\C_0$ is an irreducible component of at most two curves in $\C$. 

Summarizing, the number of $r$-rich points, with $r \ge 3$, is proportional to 
the number of incidences between the Pl\"ucker points in $\tilde{L}$ and 
the distinct curves in $\C_0$, divided by $\binom r 2$, and there is no 
contribution by the curves in $\C_1$.

\paragraph{$2$-rich points.}
The situation is different for $2$-rich points, which may arise also as incidences 
between the Pl\"ucker points of lines in $L$ and curves in $\C_1$. Handling 
them requires more care, and is done as follows. A \emph{proper} $2$-rich 
point $p$, namely a point that is incident to precisely two lines $\ell_p$ 
and $\ell'_p$ of $L$, corresponds to an incidence between the Pl\"ucker point 
of $\ell_p$ and the curve $\gamma_{\ell'_p}$ (and also between the
Pl\"ucker point of $\ell'_p$ and the curve $\gamma_{\ell_p}$). 
We count this incidence at most $\deg(\gamma_{\ell'}) = O(1)$ times, 
once for each irreducible component of the curve $\gamma_{\ell'_p}$.
It therefore suffices to count incidences between the Pl\"ucker points 
in $\tilde{L}$ and the curves in $\C_0$ (as we have just argued, this
is relevant only for curves of multiplicity at most $2$) and in $\C_1$ 
(which may have an arbitrary multiplicity). 

\medskip
Consider first incidences with curves of $\C_0$.
By projecting $T$ onto some generic plane, the number of incidences 
between the $n$ points of $\tilde{L}$ and the curves of $\C_0$ is 
the same as the number of incidences between the projected
points and the projected curves. Since $\C$ is a two-dimensional family 
of curves (Lemma~\ref{le:twodim}), so is $\C_0$. It therefore follows,
by Theorem~\ref{incPtCu} (with $s=2$), that the number of these incidences is 
$O(n^{4/3+\eps})$, for any $\eps>0$. As argued above, this gives 
us the bound $O(n^{4/3+\eps}/r^2)$ on the number of $r$-rich points,
for $r\ge 3$, thereby establishing part (a) of Theorem~\ref{th:tworich}.

This also gives us the bound $O(n^{4/3+\eps})$ for the number of $2$-rich
points that correspond to incidences formed with the curves of $\C_0$.
For the remaining $2$-rich points, which correspond to incidences between 
lines in $L$ (points of $L^*$) with 
curves of $\C_1$ (which may appear with arbitrarily large multiplicity), 
we recall that each of the curves in $\C_1$ represents one ruling of some regulus, 
and that we have assumed that no regulus contains more than $n^{1/3}$ 
lines of $L$. Hence the each curve in $\C_1$ is incident to the Pl\"ucker
points of at most $n^{1/3}$ lines in $L$, which implies that the number of 
incidences with these curves, counted with multiplicity, is at most $O(n^{4/3})$. 
Hence part (b) of the theorem also follows, and the proof is thus completed.
$\Box$

\section{Application: Distinct distances between points on an \\
algebraic curve in the plane}
\label{dd}

Let $P$ be a set of $n$ points on an irreducible algebraic curve 
$\gamma$ of constant degree in the plane, which is not a line or 
a circle. We apply the Elekes-Sharir-Guth-Katz framework~\cite{GK2}, 
and define a set $L$ of $n(n-1)$ lines in the parametric 3-space of rotations 
(rigid motions) in the plane, as $L=\{h_{a,b} \mid a\ne b \in P\}$,
where $h_{a,b}$ is the locus of all rotations that map $a$ to $b$,
which is indeed a line with a suitable parameterization; 
see \cite{ElSh,GK2} for details. Then $L$ is contained in the 
two-dimensional family of lines $\C = \{h_{x,y} \mid x,y \in \gamma\}$. 
It is easy to verify that $\C$ is irreducible and is not a $2$-flat. 
The property that $\C$ is not a 2-flat will follow from the arguments 
following Lemma~\ref{le:ref}. To see that $\C$ is irreducible, define 
a smooth morphism $\Phi: \gamma \times \gamma \to \C$ by $\Phi(x,y) = h_{x,y}$. 
If $\C$ were reducible, we could write $\C = \C_1 \cup \C_2$, and then 
$\gamma \times \gamma = \Phi^{-1}(\C_1) \cup \Phi^{-1}(\C_2)$, implying 
that $\gamma\times\gamma$ is reducible, from which it follows that $\gamma$ 
too is reducible (see, e.g., Hartshorne~\cite[Exercise 3.15(a)]{Hart83}),
contrary to assumption.

We will show below that, after pruning away some lines of $L$ (which
will not affect the asymptotic bounds derived in the analysis), the 
number of remaining lines in $L$ that are contained in a common 
plane or regulus in 3-space is $O(1)$, and thus $\C$ is 
non-degenerate with respect to $L$. Indeed, having $O(1)$ points 
in the dual space $H^*_{\ell, \ell'}$ follows by having, in primal 
space, $O(1)$ lines contained in the the plane $H_{\ell, \ell'}$.

In more detail, let $\Delta$ denote the number 
of distinct distances determined by $P$. We count the number of quadruples 
\[
\{(a,b,a',b')\in P^4 \mid |ab|=|a'b'|\},
\] 
in two different ways. First, let $N_k$ (resp., $N_{\ge k}$) denote the 
number of rotations of multiplicity exactly (resp., at least) $k$; that is, 
rotations that map exactly (resp., at least) $k$ points of $P$ to $k$ other 
points of $P$. By construction, a rotation of multiplicity at least $k$ is 
mapped to a $k$-rich point with respect to the lines of $L$. 
By Theorem~\ref{th:tworich}, replacing $\eps$ by $\eps/2$, we have 
\begin{equation} \label{ngek}
N_{\ge k} = O((n^2)^{4/3+\eps/2}/k^2) = O(n^{8/3+\eps}/k^2) , 
\end{equation}
provided that the number of lines in $L$ contained in a common 
regulus is $O(|L|^{1/3})$. We will shortly argue that this is indeed the case,
and will also show that no plane, except for a constant number of exceptional 
planes, contains more than $O(1)$ lines of $L$. 
Consider first the case of coplanar lines. A standard observation
is that one can parameterize the set of lines of the form $h_{a,b}$, 
with $a,b \in \cplx^2$, that are contained in some fixed plane, by 
rigid motions with negative determinants (i.e., motions that involve 
reflections), as follows. 

Let $\tau$ be any rigid motion of the plane with 
a negative determinant. Simple calculations show that $\tau$ must be
a reflection around some line $\ell$ in the plane, followed by a 
translation by some vector $t \in \cplx$ in the direction of $\ell$.
Such transformations are also known as \emph{glide reflections}.

We have the following lemma.
\begin {lemma} \label{le:ref}
(i) The lines of the form $h_{\xi,\tau(\xi)}$, with $\xi \in \cplx^2$, are all distinct 
and contained in a common plane, and every line in that plane is of this form, 
for a suitable point $\xi$.

\medskip
\noindent
(ii) This gives a bijection between the set of nonvertical planes in $\R^3$ 
and the set of rigid motions with negative determinants (i.e., glide reflections), 
in the sense that the plane determines the rotation $\tau$ and vice versa. 
\end {lemma}

\begin{figure}[htb]
\begin{center}
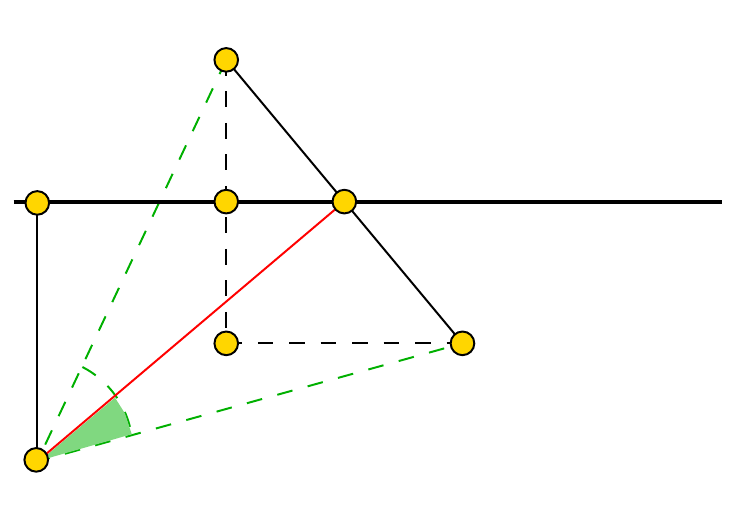
\caption{{\sf Illustrating the proof of Lemma~\ref{le:ref}: The rotation by $\theta$ about $c$
maps $a$ to $\tau(a)$, so the point $(c,\cot\frac{\theta}{2})$ lies on $h_{a,\tau(a)}$.
The similarity of $\triangle cq\mu$ and $\triangle \mu pa$ imply that
$\cot\frac{\theta}{2} = \frac{|c\mu|}{|\mu a|} = \frac{y}{t/2}$, or 
$y = \frac{2}{t} \cot\frac{\theta}{2} = \frac{2}{t} z$, which is the equation of $\pi$,
a plane that passes through $\ell$ and the angle $\beta$ that it forms with the $xy$-plane
satisfies $\tan\beta = \frac{t}{2}$.}} 
\label{fig:glide}
\end{center}
\end{figure}

\noindent
{\bf Sketch of proof.}
The plane $\pi$ in part (i) is obtained as follows. Its intersection with the 
$xy$-plane is the line $\ell$, and the angle $\beta$ that it forms with the 
$xy$-plane is given by $\tan\beta = \frac{t}{2}$. See Figure~\ref{fig:glide}
for an illustration of this claim. This also gives a recipe to reconstruct $\tau$ from $\pi$.
$\Box$

Assume that the lines $h_{a_i, b_i}$, for $i=1,\ldots, r$, are contained 
in some non-vertical plane. Then, by Lemma~\ref{le:ref}, there is a glide reflection
$\tau$, such that $h_{a_i, b_i} = h_{a_i, \tau(a_i)}$ (this latter property 
follows because $h_{a,b} \ne h_{a',b'}$ when $(a,b) \ne (a',b')$), and thus 
$b_i = \tau(a_i)$, for $i=1,\ldots,r$. However, this implies that 
$a_1,\ldots, a_r \in \gamma \cap \tau^{-1}(\gamma)$, which implies that 
$\gamma = \tau(\gamma)$ (assuming that $r$ is larger than the degree of 
$\gamma$). By Pach and de Zeeuw \cite[Lemma 2.5]{PdZ}, 
this can happen for at most $4\deg(\gamma)=O(1)$ rotations, unless $\gamma$ 
is a line orthogonal to or coinciding with $\ell$ (assuming $\gamma$ to be 
irreducible), or $\gamma$ is a circle and $\ell$ goes through its center. 
Since both these possibilities have been excluded, we are left with
at most $O(1)$ glide reflections $\tau$ (that is, planes in 3-space) that satisfy
$\gamma = \tau(\gamma)$. Fix one of these planes $\pi$, with a corresponding
glide reflection $\tau$, and note that the number of lines $h_{a,\tau(a)}$
is at most $n$. Hence there are at most $O(n^2)$ rich points (that is, $2$-rich points)
formed solely by lines in $\pi$. Any other rich point occurs at an intersection 
of a line of $L$ (not contained in $\pi$) with $\pi$, and the number of these points
is $O(n^2)$. To recap, the $O(1)$ exceptional rotations contribute at most
$O(n^2)$ $r$-rich points, for any $r\ge 2$, well below our bound. We remove
the lines lying in these $O(1)$ planes, and, among the remaining lines,
no plane contains more than $O(1)$ lines.

This argument also implies that $\C$ cannot be a 2-flat. Indeed, 2-flats
within the Klein quadric $Q$ come either from lines in $\reals^3$ that are
contained in a common plane or from lines that pass through a common point
(see Rudnev and Selig~\cite{RS} for details). The former situation is
impossible because, as we have argued, it is impossible for all (or just 
for too many of) the lines $h_{x,y}$, for $x,y\in\gamma$, to be contained 
in a common plane in the primal 3-space. The latter situation is also easy,
as the point in $\reals^3$ common to the lines $h_{x,y}$ is a rotation that 
maps every point $x\in\gamma$ to every other point $y\in\gamma$, which is 
clearly impossible.

The situation is simpler for lines in a common regulus. Assume that the 
lines $h_{a_i, b_i}$, for $i=1,\ldots, r$, are contained in some regulus 
$\sigma$. In general, any regulus in $\cplx^3$ is equivalent to the 
hyperbolic paraboloid $z = xy$ by a suitable change of coordinates, 
so we may assume that $\sigma$ is $z = xy$. The lines contained in 
$\sigma$ are either of the form $x = a, z = ay$ or of the form 
$y=a, z=ax$, for $a\in \cplx$. We may assume, without loss of generality,
that our lines are of the first kind, and then it is easy to verify that 
these lines can be expressed as $h_{(a-\frac 1 a, 0), (a + \frac 1 a, 0)}$, 
for $a\in\cplx$. Hence, if there are at least $k$ values $a$, where $k$
is the degree of $\gamma$, such that both points $(a-\frac 1 a, 0)$ and 
$(a + \frac 1 a, 0)$ belong to $\gamma$, then, since $\gamma$ is 
irreducible and of constant degree, it must coincide with the $x$-axis, 
contradicting our assumption that $\gamma$ is not a line.

We have thus shown that no plane or regulus contains more than a constant
number of lines of $L$, except for at most $O(1)$ special planes, whose
effect on the asserted bound is negligible, and which we ignore by removing
all the lines contained in these planes. 
%
Hence, arguing as in \cite{ElSh,GK2} and using (\ref{ngek}), 
the number of quadruples is at most
\begin{equation}
\sum_{k=2}^n \binom k 2 N_k \le \sum_{k=2}^n k N_{\ge k} 
= O\left( \sum_{k=2}^{n} \frac{n^{8/3+\eps}}{k} \right) = O(n^{8/3+\eps}\log n) = O(n^{8/3+2\eps}).
\end{equation}
On the other hand, by Elekes's analysis, which is based on the Cauchy-Schwarz 
inequality (see, e.g., Guth and Katz~\cite{ElSh,GK2}), the number of quadruples 
is also $\Omega(n^4/\Delta)$, implying that the number of distinct distances 
satisfies $\Delta=\Omega(n^{4/3 -2\eps})$. This result was obtained earlier 
in~\cite{PdZ}, without the $\eps$-loss in the exponent, but the proof here is 
much simpler, and we hope that it will find similar applications of this kind.

\section{Incidences between points and lines in a two-dimensional family of lines} \label{sec:inc2d}

The main result of this section is the following theorem.

\begin{theorem}
\label{th:inctw}
Let $P$ be a set of $m$ points in $\R^3$, and let $L$ be a set of $n$ 
lines in $\R^3$, such that the set $\tilde{L}$ of the Pl\"ucker images of the
lines of $L$ is contained in some two-dimensional, non-planar, irreducible 
variety $T$ of constant degree, as in Section~\ref{rich-2d}. Then $I(P,L)=O(m+n)$.
\end{theorem}

\noindent{\bf Proof.}
As observed above, for a point $p \in \cplx^3$, the set $S_p$ of the Pl\"ucker 
images of lines that are incident to $p$ form a $2$-flat in Pl\"ucker space, 
and that 2-flat is contained in $Q$. We first observe that even if $T$ were a 
$2$-flat in $Q$, and, for some $p \in \R^3$, we had $T = S_p$, then 
all the lines with Pl\"ucker images in $T$ would be incident to $p$, 
and thus for any point $q \ne p$, the number of incidences between $q$ 
and the lines with Pl\"ucker images in $T$ would be at most one 
(there is only one line in $\cplx^3$ that is incident to $p$ and $q$), 
so the number of incidences would be $O(m+n)$ in this case too. 
Thus, we may assume that for every $p\in \cplx^3$, $T \ne S_p$. 
As $T$ is two-dimensional, B\'ezout's theorem~\cite{Fu84} implies 
that for every $p\in \cplx^3$, the intersection $S_p \cap T$ is a 
union of a constant number of curves of constant degree and a discrete 
set of a constant number of points, in Pl\"ucker space. 

Put 
\[
V := \{p \in \cplx^3 \mid S_p \cap T \text{ is a curve}\}.
\] 
Similar to~\cite[Theorem 2.16]{SS4d}, one can define a polynomial of constant degree, 
via multivariate resultants, whose vanishing at a point $p\in\cplx^3$ is equivalent 
to $S_p \cap T$ being one-dimensional (in $Q$). Hence, $V$ is a complex algebraic 
variety, and, as $T$ is of constant degree, so is $V$. We assume, as we may, that $V$ is
irreducible (by treating each of the $O(1)$ irreducible components of $V$ separately).

Since $T$ is two-dimensional (in $Q$), $V \ne \cplx^3$. Indeed, write
$T = \bigcup_{p \in \cplx^3} (S_p \cap T)$ (the equality holds because
$\bigcup_{p \in \cplx^3} S_p$ is the entire Klein quadric $Q$).
Then $V=\cplx^3$ would imply that 
$\dim(T) \ge 3$, by counting one dimension, namely $\dim(S_p \cap T)$, for every 
point $p \in \cplx^3$, and then subtracting one dimension, as $S_p \cap T$ is 
counted for every point in $S_p \cap T$, which is a one-dimensional curve in $Q$, 
thus contradicting the fact that $T$ is $2$-dimensional.

We next argue that $V$ cannot be two-dimensional (recall that $V\subset\cplx^3$). 
Assume to the contrary that $V$ is two-dimensional. By definition, for each 
$p\in V$, we have that $T$ contains a one-dimensional set of the Pl\"ucker images
of lines that are incident to $p$ (namely, those images forming the curve $S_p\cap T$ 
in the Pl\"ucker space).

\begin{lemma}
\label{le:int}
If $V$ is a 2-flat in the primal $\cplx^3$, then $T$ is a $2$-flat in $Q$. 
\end{lemma}

\noindent{\bf Proof.}
For a point $p \in \cplx^3$, let
$X_p:= \{\tilde{\ell} \in S_p \cap T \mid \ell \subset V \} = S_p \cap T \cap V^*$, where
$V^* = \{ \tilde{\ell} \mid \ell\subset V\}$.
By B\'ezout's Theorem, $S_p \cap T$ is a union of $O(1)$ constant-degree curves 
contained in $S_p$, and $O(1)$ Pl\"ucker images of lines (note that $S_p \cap T$ 
can also be a discrete finite set). Let 
\[
U := \{p \in V \mid X_p \text{ is a curve in Pl\"ucker space}\}.
\] 
We claim that $U$ is a subvariety of $V$. To see this,
define $Y = \{(p, \ell) \in V \times (T \cap V^*) \mid \ell \in S_p\}$. 
Clearly $Y$ is an algebraic variety whose fiber over point $p\in V$ of the projection onto $V$ is $Y_p = T \cap V^*\cap S_p$, and by the 
Theorem of the Fibers~\cite[Corollary 11.13]{Har} 
(see also Theorem 6.2 in~\cite[Theorem 6.2]{SS4dv}), we deduce that 
$\{p\in V \mid \dim(Y_p)\ge 1\}$ is an algebraic variety, and this is precisely $U$.


If $U$ is two-dimensional, then $V=U$, implying that for every point $p$ in 
the 2-flat $V$, there are infinitely many lines $\ell$ that are incident to 
$p$ and contained in $V$ and that $\ell^*. \in T$. This implies that 
$T = V^*$, and that $T$ is a $2$-flat. 

Otherwise, $U$ is at most one-dimensional, and for every 
$p \in V\setminus U$, there is a set $K$, consisting of $O(1)$ Pl\"ucker images
$\tilde{\ell}$ in $S_p \cap T$ of lines that are contained in $V$ (in the primal 
3-space). For a pair $p \ne q \in V \setminus U$, $S_p$ and $S_q$ are both 
$1$-dimensional (in $Q$), and, outside the set $K$ of constantly many Pl\"ucker 
images of lines, $S_p$ and $S_q$ are disjoint. Thus, $\dim(T) = 3$, as it is 
parameterized by three parameters, two for the point $q \in V$, and one for the
$1$-dimensional curve $S_p$, contradicting the assumption that $\dim(T)=2$.
$\Box$

\medskip

We formalize the argument from Lemma~\ref{le:int}, and generalize it to the 
general two-dimensional case, as follows. For a point $p \in \cplx^3$, the 
set $X_p:= \{\tilde{\ell} \in S_p \cap T \mid \ell \subset V\} = S_p \cap T \cap V^*$ is a union 
of $O(1)$ constant-degree curves contained in $S_p$, and $O(1)$ additional 
Pl\"ucker images of lines. 

As observed in \cite[Appendix]{SS4dv}, a surface is ruled by lines if every 
point in a Zariski-open dense subset on the surface is incident to a line 
that is fully contained in it. The same is true for infinitely ruled surfaces, 
namely, it suffices to have a Zariski-open dense subset on the surface, each 
of whose points is incident to infinitely many lines that are fully contained in 
the surface. Therefore, if the set $U$, defined in the proof of Lemma~\ref{le:int},
is two-dimensional, it follows that $V$ is infinitely ruled by lines, and by 
\cite[Theorem 3.11]{SS4d}, it follows that $V$ must be a 2-flat (in the primal 
3-space), and Lemma~\ref{le:int} implies that $T$ is a $2$-flat in $Q$, 
contradicting our assumption. Otherwise, $U$ is at most one-dimensional, and 
for every $p \in V\setminus U$, there are only $O(1)$ Pl\"ucker images $\tilde{\ell}$
in $S_p \cap T$ of lines $\ell$ that are contained in $V$. Hence, for any pair 
$p \ne q \in V \setminus U$, $S_p$ and $S_q$ are both $1$-dimensional and, 
outside a possible set $K$ of $O(1)$ Pl\"ucker images of lines that are 
contained in $V$, $S_p$ and $S_q$ are disjoint. It thus follows that $\dim(T) = 3$, 
as it is parameterized by three parameters, two for the point $q \in V$, and one 
for the $1$-dimensional curve $S_p$, contradicting the assumption that $\dim(T)=2$. 
Thus $V$ is a curve, contradicting our assumption that $V$ is two-dimensional.

We have thus argued that $V$ cannot be two-dimensional, so $V$ must be one-dimensional.
By definition of $V$, every point $p\in P \setminus V$ is incident to at most $O(1)$
lines of $L$, for a total of $O(m)$ incidences.
Thus, we may assume that all the points of $P$ are contained in the curve $V$.
For each line $\ell \in L$, if $\ell$ is not contained in $V$ it contributes at most 
$O(\deg V) = O(1)$ incidences with $P$. Thus, we get a total of $O(n)$ incidences, except for 
at most $O(\deg V) = O(1)$ lines that are contained in $V$, for a total of $O(m)$ additional incidences.

This completes the proof of the theorem.
$\Box$

\medskip
The following corollary is an immediate consequence of the theorem.
\begin{corollary}
\label{co:easy}
Let $T$ be a two-dimensional, non-planar, irreducible subvariety, of 
constant degree, of the Klein quadric $Q$. Then, there exists a constant 
$r_0 = r_0(\deg(T))$ so that, if $L$ is a set of $n$ lines in $\R^3$ 
whose Pl\"ucker images are points in $T$ then, for $r \ge r_0$, the set 
$P_{\ge r}(L)$ of $r$-rich points determined by $L$ satisfies $|P_{\ge r}(L)| = O(n/r)$.
\end{corollary}

\section{Incidences between points and lines in a three-dimensional family of lines}
\label{sec:inc:3d}

In this section we prove the following result.
\begin{theorem} \label{incmain}
The number of incidences between $m$ points in $\R^3$ and $n$ lines in $\R^3$
whose Pl\"ucker images are contained in an irreducible nonlinear constant-degree 
three-dimensional variety $T$ is 
\[
O\left(m^{3/5}n^{3/5} + (m^{11/15}n^{2/5} + m^{1/3}n^{2/3})s^{1/3} + m + n \right) ,
\]
provided that no plane contains more than $s$ of the points. 
If $s = O(\min\{n^{3/5}/m^{2/5}, m^{1/2}\})$, the bound becomes $O(m^{3/5}n^{3/5} + m + n)$.
\end{theorem}

\noindent{\bf Proof.}
Recall that $S_p$ is the 2-flat in Pl\"ucker space that consists of the 
images of all lines passing through a point $p\in\R^3$, and define 
\[
W := \{p \in \cplx^3 \mid S_p \cap T \text{ is two-dimensional}\}.
\] 
\begin{lemma}
\label{le:td}
$W$ is an algebraic variety of dimension at most $2$ and of constant degree.
\end{lemma}

\noindent{\bf Proof.}
Since $S_p$ is a 2-flat, $S_p \cap T$ is two-dimensional 
if and only if $S_p \subset T$. Similarly to the Fano
variety of lines, the Grassmannian manifold of $2$-flats contained in a
constant-degree variety is an algebraic variety of constant 
degree \cite{GrHa}, so $W$ is algebraic of constant degree. To bound the dimension 
of $W$, we repeat the proof of \cite[Theorem 2.3(a)]{SS4dv}, 
which proceeds by counting the dimensions of the fibers that arise in
the problem. Here we omit the details and give the high-level idea. 
Assume to the contrary that $W$ is three-dimensional, i.e., 
$W = \cplx^3$, so for every point $p\in\cplx^3$, the 
$2$-flat $S_p$ is contained in $T$. Omitting details, we note 
that each $p \in \cplx^3$ contributes a two-dimensional set ($\dim(S_p)=2$), 
but then every line is counted by the infinitely many points incident to it.
A standard dimension counting argument then implies that
$\dim(F(T)) \ge 4$, where $F(T)$ is the Fano variety 
of lines contained in $T$. By \cite[Theorem 3.11]{SS4d}, this implies
that $T$ has to be a $3$-flat, contradicting our assumption.
$\Box$

\medskip
We first treat incidences with points $p \in P \cap W$. We decompose $W$ 
into its $O(1)$ irreducible components, and treat each component separately. 
If a component $W_0$ of $W$ is not a $2$-flat then, by~\cite[Corollary 1.4]{SS4dv}, 
the number of incidences between points contained in $W_0$ and lines in $L$ is $O(m+n)$. 
If $W_0$ is a $2$-flat, we invoke the Szemer\'edi-Trotter bound in 
Theorem~\ref{th:ST}, and get the bound $O(s^{2/3}n^{2/3} + s + n)$,
using our assumption that no 2-flat contains more than $s$ points of $P$.
This in turn can be upper bounded by $O(s^{1/3}m^{1/3}n^{2/3} + m + n)$, 
which is subsumed by the bound asserted by the theorem.

Next, we treat incidences involving points in $\R^3 \setminus W$, i.e., 
points that are incident to a one-dimensional family of lines with Pl\"ucker images in $T$. 
In this case we use duality, replacing each point $p$ in $\cplx^3$ with 
the one-dimensional curve $\gamma_p$ of the Pl\"ucker images $\tilde{\ell}$
of lines $\ell$ that are incident to $p$ and have Pl\"ucker images in $T$.
This yields a family of $m$ constant-degree curves that is a family of 
pseudo-lines. (Two such curves $\gamma_p$ and $\gamma_q$ intersect in at 
most one point, representing the (unique) line connecting $p$ and $q$, if 
its Pl\"ucker image lies in $T$.) We replace each of the $n$ lines 
in $L$ by its Pl\"ucker image, and obtain an incidence problem between $n$ points 
and $m$ pseudo-lines within the variety $T$, a three-dimensional subset of the 
four-dimensional Klein quadric $Q$. Using a generic projection of $T$ onto $\R^3$
(in which all projected points are distinct and no pair of projected curves overlap),
the analysis then proceeds by invoking Zahl~\cite[Lemma 4.1]{Za17},
which extends the Guth--Katz incidence bound from incidences with lines 
to incidences with pseudo-lines. Specifically, Zahl shows that the
number of incidences between $n$ points and $m$ pseudo-lines in $\R^3$, 
assuming that these pseudo-lines are constant-degree algebraic curves, is
$O(n^{1/2}m^{3/4}+n^{2/3}m^{1/3}\xi^{1/3} + m + n)$, where $\xi$ is an upper 
bound on the number of pseudo-lines that are contained in any common 
two-dimensional surface contained in $T$ that is infinitely ruled by curves 
from the infinite family from which our pseudo-lines are taken.

As argued in Guth and Zahl~\cite{GZ}, any such surface must be of degree
at most $100E^2$, where $E$ is the degree of the pseudo-lines $\gamma_p$, 
so it is sometimes convenient, especially when no simple characterization of
such infinitely ruled surfaces is known, to impose the stronger assumption
that no surface of degree at most $100E^2$ contains more than $\xi$ pseudo-lines.
In general, without having a good characterization of the infinitely-ruled
surfaces, this assumption is too restrictive, and difficult to verify.
One of the main technical contribution of the analysis in this section 
is to exploit the dual nature of the present setup, and replace this 
assumption by the simpler and more natural assumption that, in the 
original primal 3-space, \emph{no plane contains more than $s$ points 
of $P$}, allowing us to replace $\xi$ by $s$.

To see this we argue as follows.
Let $S$ be a surface contained in $Q$ that is infinitely ruled by curves from $\C_0$, 
and define the variety\footnote{%
  It is easy to verify that $V_S$ is defined by polynomial equations and is thus an algebraic variety.} 
$V_S := \{p \in \cplx^3 \mid \text{an irreducible component of }\; c_p \subseteq S\}$. We have
\begin{lemma}
$V_S$ is two-dimensional.
\end{lemma}

\noindent{\bf Proof.}
We assume, as we may, that $S\subset T$, and note that then
$V_S = \{p \in \cplx^3 \mid \dim(S_p \cap S)=1\}$, and that $V_S$ parameterizes curves
$c_p$ that are contained in $S$. By using resultant theory (see, e.g.,
\cite[Section 2.4]{SS4d}), it is also easy to verify that $V_S$ is an algebraic variety.

Since $S$ is infinitely ruled by curves in $\C_0$, it follows that for
every $\tilde{\ell} \in S$, there exists an infinite family of curves in $\C_0$
that are contained in $S$ and incident to $\tilde{\ell}$. Since every curve in
$\C_0$ is of the form $S_p \cap S$ with $p \in \ell$, it follows that 
each point $\tilde{\ell}\in S$ is incident to a $1$-dimensional family of curves 
of $\C_0$ (which is the family corresponding to points in $\ell$).
As $S_p \cap S_{q}$ is a point, it follows that pairs of curves in $\C_0$
intersect in at most one point.
By dimension counting, it follows that $V_S$ is two-dimensional.
$\Box$

\medskip

By definition of $S$, for every point $\tilde{\ell} \in S$, there are infinitely many curves 
of the form $c_p$, for $p \in V_S$, that are contained in $S$ and incident to $\tilde{\ell}$. 
For each such curve $c_p$, in the primal space $\cplx^3$, the corresponding point $p$ 
is incident to the line $\ell$. Therefore, there are infinitely many points of $\ell$ that 
are contained in $V_S$, and as $V_S$ is an algebraic variety, it follows that $\ell \subset V_S$.

Therefore, we have a two-dimensional variety $V_S \subset \cplx^3$ containing 
lines that are parameterized by a \emph{two-dimensional} variety $S \subset T$, 
implying that $V_S$ is a plane. (Indeed, the Fano variety of lines contained in 
a two-dimensional irreducible variety that is not a plane is at most one-dimensional.)

In summary, the lines in $\cplx^3$ whose Pl\"ucker coordinates lie in $S$ are 
contained in a plane $V_S$. As observed above (and argued in Rudnev~\cite{Rud}), 
the set of lines that are contained in a plane in $\cplx^3$ are dual to a 2-flat 
in the Klein quadric $Q$, implying that $S$ is a plane.

For a point $p \in V_S$, we have $S_p \cap T = c_p \subset S$, implying that 
$S_p \cap T = S_p \cap S$, and the right-hand side is an intersection of 
two 2-flats. Except for at most one point $p_0$, we get that for every (other) 
$p \in V_S$, the curve $S_p \cap S$ is a \emph{line} contained in $S$. 

Since the set of lines fully contained in $S$ is two-dimensional, it follows\footnote{%
  The set of lines of the form $S_p \cap S$ is a two-dimensional variety 
  contained in the set of lines contained in $S$, which is an irreducible 
  two-dimensional variety. Moreover, for two distinct points $p,q \in V_S$, 
  the intersection $S_p \cap S_q$ is equal to the point in 
  Pl\"ucker space that represent the unique line that is incident to $p$ and $q$. 
  Therefore, the lines $S_p \cap S$, for $p\in V_S$, are distinct. 
  Thus the two varieties must be equal.} 
that all the lines fully contained in $S$ are of the form $S_p \cap S$. 
By assumption, the number of points in $P$ that are contained in a common plane 
is at most $s$, and duality implies that the number of lines of the form 
$S_p \cap S$ fully contained in $S$ is at most $s$. 

We thus obtain the incidence bound
\[
I(P,L) = O(n^{1/2}m^{3/4}+n^{2/3}m^{1/3}s^{1/3} + m +n) .
\]
Since $n^{1/2}m^{3/4} \le m^{3/2}$ when $n \le m^{3/2}$, and $n^{1/2}m^{3/4} \le n$ 
otherwise, we get the following bootstrapping bound
\begin{equation}
    \label{eq:bootb}
    I(P,L) = O(m^{3/2} + n^{2/3}m^{1/3}s^{1/3} + n) .
\end{equation} 
The analysis then proceeds by ``starting over'' in primal space, 
i.e., by constructing a partitioning polynomial $g$ of degree $O(D)$,
for a suitable value of $D$, to be fixed shortly, using the techniques 
in~\cite{Guth,GK2}, so that each connected component (cell) $\tau$ of 
$R^3 \setminus Z(g)$ contains at most $m/D^3$ points of $P$ and is crossed 
by at most $n/D^2$ lines of $L$ (but any number of points and lines can
be contained in the zero set $Z(g)$). 

\paragraph{Incidences within the cells.}
We first bound the number of incidences within the partition cells.
We apply the bootstrapping bound in (\ref{eq:bootb}) to each cell 
$\tau$ and sum the bound over all components, to obtain the bound 
 \begin{multline*}
 O\left( D^3 ((m/D^3)^{3/2} + (n/D^2)^{2/3}(m/D^3)^{1/3}s^{1/3} + n/D^2)\right) \\
 = O\left( \frac{m^{3/2}}{D^{3/2}} + n^{2/3}m^{1/3}D^{2/3}s^{1/3} + nD\right) .
 \end{multline*}
To balance the first and last terms, we choose $D = m^{3/5}/n^{2/5}$.
For this to make sense, we require that $1\le D\le \min\{m^{1/3}, n^{1/2}\}$, 
or, equivalently, that $n\le m^{3/2}$ and $m\le n^{3/2}$. 
When the first inequality does not hold, we do not use any 
partitioning and just apply (\ref{eq:bootb}) to obtain the bound
$I(P,L) = O(n^{2/3}m^{1/3}s^{1/3} + n)$. When the second inequality 
does not hold, we choose $D = an^{1/2}$, for a suitable constant $a$,
which satisfies the inequalities $1\le D\le \min\{m^{1/3}, n^{1/2}\}$. 
In fact, we can construct a polynomial
$g$ of this degree so that all the lines of $L$ are fully contained in
$Z(g)$ (see, e.g.,~\cite{KSS}),
and we may therefore assume that all the points of $P$ are also contained in $Z(g)$, as the other 
points contribute no incidences.  That is, in this case there are no incidences within the cells.

In the middle range, our choice of $D$ yields the bound 
$O(m^{3/5}n^{3/5} + m^{11/15}n^{2/5}s^{1/3})$. Combining all the bounds,
the number of incidences within the partition cells is
\begin{equation} \label{inc:incells}
O\left(m^{3/5}n^{3/5} + (m^{11/15}n^{2/5} + m^{1/3}n^{2/3})s^{1/3} + n \right) .
\end{equation}

\paragraph{Incidences on the zero set.}
Consider next incidences involving points that lie on $Z(g)$. 
A line $\ell$ that is not fully contained in $Z(g)$ crosses it in
at most $O(D)$ points, for an overall $O(nD)$ bound, which is subsumed 
by the bound (\ref{inc:incells}) for incidences within the cells.
It therefore remains to bound the number of incidences between the points 
of $P$ on $Z(g)$ and the lines that are fully contained in $Z(g)$. 

We handle each irreducible component of $Z(g)$ separately. 
For non-planar components, Theorem \ref{th:ss4dv}, combined with
H\"older's inequality (for summing up the bounds over the irreducible 
components) implies that the number of incidences between points and 
lines contained in $Z(g)$, but not in any planar component of $Z(g)$, is 
\[
O(m^{1/2}n^{1/2}D^{1/2} + m^{1/3}D^{4/3}s^{1/3} + m + n) .
\] 
In the middle range $n^{2/3} \le m \le n^{3/2}$, the choice of
$D = m^{3/5}/n^{2/5}$ is easily seen to yield the desired bound
$O(m^{3/5}n^{3/5} + m + n)$ (as the second term in the bound is dominated
by the first term for this range of $m$). The case $m < n^{2/3}$ has already been 
handled, by a single application of (\ref{eq:bootb}), which
yields the bound $O(n^{2/3}m^{1/3}s^{1/3} + n)$.
When $m > n^{3/2}$, the choice of $D = an^{1/2}$, as made above, 
yields the bound $O(n^{2/3}m^{1/3}s^{1/3} + m)$.

For the planar components, we use the standard technique of assigning 
each point and line to the first planar component that contains it
(according to some arbitrary enumeration of the components). The
number of incidences between points and lines assigned to different 
components is $O(nD) = O(m^{3/5}n^{3/5}+m+n)$ (the right-hand side
does indeed bound the left-hand side for each of the sub-ranges). 
For incidences between
points and lines assigned to the same planar component, we apply 
the Szemer\'edi-Trotter bound (Theorem \ref{th:ST}) to each component 
and sum the resulting bounds over the components. The assumption that 
each plane contains at most $s$ points, combined with H\"older's 
inequality, yields the bound $O(m^{1/3}n^{2/3}s^{1/3} + m + n)$.

That is, the number of incidences with points on $Z(g)$ is bounded by
\begin{equation} \label{inc:onzg}
O\left( m^{3/5}n^{3/5} + m^{1/3}n^{2/3}s^{1/3} + m + n \right) .
\end{equation}
Combining with the bound (\ref{inc:incells}) for incidences within the cells,
we get the overall bound
\[
I(P,L) = O\left(m^{3/5}n^{3/5} + (m^{11/15}n^{2/5} + m^{1/3}n^{2/3})s^{1/3} + m + n \right) ,
\]
thereby completing the proof of the theorem.
$\Box$

\medskip

\noindent{\bf Remark.} 
An interesting open challenge in incidence geometry is to sharpen
the Guth-Katz bound \cite{GK2} when the number of lines in any common plane
is at most some constant. When the lines in $L$ are contained, as points
in Pl\"ucker space, in an irreducible nonlinear constant-degree 
three-dimensional variety $T$ then, while we cannot deduce that the 
number of lines contained in a common plane is constant, we can 
nevertheless deduce the following useful property. For any plane 
$\Pi \subset \cplx^3$, $T \cap \Pi^*$ (recall that $\Pi^*$ is the 
2-flat dual to $\Pi$, consisting of all the points dual to lines 
that are contained in $\Pi$) is a constant-degree curve, and thus, 
except for $O(1)$ points, every point in $\Pi$ is incident to $O(1)$ 
lines in $T$, implying that the number of incidences in a common 
plane is $O(m_{\Pi}+n_{\Pi})$, where $m_{\Pi}$ ($n_{\Pi}$) is the 
number of points (lines) contained in $\Pi$. Such a linear bound
on the number of incidences within a plane is a key property for 
deriving improved incidence bounds, as demonstrated in this work.
For Theorem \ref{incmain}, we also added the condition that 
$m_{\Pi} \le s$, for every plane $\Pi$, to further improve the bound.

\section{Application: Incidences between points and lines on a quadric in four dimensions}

Solomon and Zhang~\cite{SoZ} give a configuration of $m$ points and $n$ lines in 
a quadratic hypersurface in $\R^4$, having $\Omega(m^{2/3}n^{1/2}+m+n)$ incidences. 
The following theorem follows as a corollary from the previous section.
\begin{theorem}
\label{th:incfour}
Let $P$ be a set of $m$ points and $L$ a set of $n$ lines contained in a quadratic 
hypersurface $S\subset \C^4$ such that no $2$-flat contains more than $s = O(n^{3/5}/m^{2/5})$ 
of the points of $P$. Then $I(P,L)=O(m^{3/5}n^{3/5} + m+ n)$.
\end{theorem}

\medskip
\noindent{\bf Remark.}
When $m = O(n^{3/2})$, the lower bound $\Omega(m^{2/3}n^{1/2}+m+n)$ obtained in 
\cite{SoZ} is (asymptotically) smaller than the upper bound $O(m^{3/5}n^{3/5}+m+n)$ 
asserted in Theorem~\ref{th:incfour}. Closing this gap remains a challenging open problem.

\paragraph{Acknowledgements.}
We deeply thank Larry Guth for helpful interaction on this work.
We also thank an anonymous referee for his/her many careful and constructive comments.
Having addressed them has greatly improved the presentation of this work.

\end{document}